\documentclass{amsart}
\usepackage{palatino,hhline, bbm}
\usepackage{arydshln} 
\usepackage{amsthm, amsmath}
\usepackage{amssymb} 
\usepackage{amscd}
\usepackage{mathrsfs}
\renewcommand{\mathcal}{\mathscr}
\usepackage[hmargin=3cm, vmargin=1.8cm]{geometry}
\usepackage[breaklinks,bookmarksopen,bookmarksnumbered]{hyperref}
\usepackage[all]{xy}

\theoremstyle{plain}
\newtheorem*{thm}{Theorem}
\newtheorem*{lem}{Lemma}

\theoremstyle{remark}

\newcommand\rond{\kern 1pt{\scriptstyle\circ}\kern 1pt}

\def\lr#1{\langle {#1} \rangle}

\newcommand\im{\operatorname{Im}}
\newcommand\Ker{\operatorname{Ker}}
\newcommand\Coker{\operatorname{Coker}}
\newcommand\Pic{\operatorname{Pic}}
\newcommand\Div{\operatorname{Div}}

\newcommand\Z{\mathbb{Z}}
\newcommand\Q{\mathbb{Q}}
\newcommand\R{\mathbb{R}}
\newcommand\C{\mathbb{C}}
\renewcommand\P{\mathbb{P}}

\renewcommand\O{\mathcal{O}}

\newcommand\mono{\lhook\joinrel\mathrel{\longrightarrow}}
\newcommand\iso{\vbox{\hbox to .8cm{\hfill{$\scriptstyle\sim$}\hfill}
\nointerlineskip\hbox to .8cm{{\hfill$\longrightarrow $\hfill}} }}
\newcommand\bir{\vbox{\hbox to .8cm{\hfill{$\scriptstyle\sim$}\hfill}
\nointerlineskip\hbox to .8cm{{\hfill$\dasharrow $\hfill}} }}
\newcommand\abs[1]{\lvert {#1}\rvert}

\begin{document}
\title{Even sets of nodes and Gauss genus theory}
\author[Arnaud Beauville]{Arnaud Beauville}
\address{Universit\'e C\^ote d'Azur\\
CNRS -- Laboratoire J.-A. Dieudonn\'e\\
Parc Valrose\\
F-06108 Nice cedex 2, France}
\email{arnaud.beauville@unice.fr}
 
\begin{abstract}
We observe that a lemma used in the study of even sets of nodes on surfaces applies almost verbatim to prove a celebrated formula of Gauss on the $2$-torsion of the class group of a quadratic field.
\end{abstract}
\maketitle 
\section{Introduction}
It is a great pleasure for me to dedicate this paper to Herb Clemens. I met Herb in Chile in 1972, and this played a decisive role in my mathematical orientation.

 I visited Herb for one month in Salt Lake City in 1979. I was interested in surfaces at that time; I~gave in particular a talk about the maximum number of double points that a quintic surface can have. A crucial ingredient was a lemma expressing the 2-torsion of the Picard group of a double cover in terms of the branch locus.
 

Gauss genus theory deals with binary quadratic forms; one of its main consequences is the precise description of the 2-torsion of the ideal class group of a quadratic field in terms of the discriminant.
It is only quite recently that I realized that the statement is very close to that of my lemma -- in fact, the   lemma applies almost verbatim to give a proof of Gauss' theorem. This is what I would like to explain in this note.

\emph{Disclaimer}: That proof is not new, and can be found (with a somewhat different language) in many number theory books. So this note is essentially expository, with the aim of highlighting  an unexpected connection between a classical question  in algebraic geometry and a basic result in number theory.

\bigskip	
\section{Gauss' theorem}
Let $K$ be a number field, and $\O$ its ring of integers. The \emph{ideal class group} $\operatorname{Cl}(K)$ is the Picard group of $\operatorname{Spec}(\O) $;  we will denote it simply  by $\Pic(\O)$. It is the quotient of the group of fractional ideals $\operatorname{Div}(\O) $ (the free abelian group with basis the nonzero prime  ideals of $\O$) by the subgroup of principal ideals. Thus we have
an exact sequence
\[1\rightarrow \O^*\rightarrow K^*\rightarrow \operatorname{Div}(\O)\rightarrow \Pic(\O)\rightarrow 0 \,.\]

An important variant of this definition takes into account the places  at infinity of $K$. Let $K^*_+$  be the subgroup of elements  $\alpha \in K^*$  which are \emph{totally positive}, that is, such that $\sigma (\alpha )>0$ for all embeddings $K\mono \R$; put $\O^*_+:=\O\cap K^*_+$. Replacing $K^*$ by $K^*_+$ in the above definition we get the narrow class group $\operatorname{Cl}_+(K)=\Pic_{+}(\O) $, which fits into an exact sequence
\[1\rightarrow \O^*_+\rightarrow K^*_+\rightarrow \operatorname{Div}(\O)\rightarrow \Pic_+(\O)\rightarrow 0 \,.\]
As a consequence of his deep study of binary quadratic forms, Gauss obtained the following:
\begin{thm}[Gauss]
Let $K$ be a quadratic extension of $\Q$, ramified at $r$ primes. Then\footnote{ For an abelian group $G$, we denote by $G[2]$ its  2-torsion part. } $\operatorname{Cl}_+(K)[2]=(\Z/2)^{r-1} $.
\end{thm}
$\bigl($If $K=\Q(\sqrt{d})$ and $d=\pm \,p_1\ldots p_s$, the ramified primes are 
$p_1,\ldots ,p_s$, plus $2$ if $d\equiv 3\pmod{4}$.$\bigr)$

\smallskip	
From the   Remark below one gets  $\operatorname{Cl}(K)[2]\cong (\Z/2)^{r-1}$ or $(\Z/2)^{r-2}$. Gauss' theorem is quite remarkable: such a simple expression does not exist for the $p$-torsion of $\operatorname{Cl}(\Q(\sqrt{d})) $ for $p>2$, or the 2-torsion of $\operatorname{Cl}(K) $ for $[K:\Q]>2$. 

\medskip	
\noindent\emph{Remark}$.-$ When $d<0$, we have $K^*_+=K^*$ and $\operatorname{Cl}_+(K)=\operatorname{Cl}(K)$. When $d>0$,
there is a surjective map $\operatorname{Cl}_+(K)\rightarrow \operatorname{Cl}(K)$; its kernel is  $0$ or $\Z/2$, according whether $\O$ contains elements of norm $-1$ or not. Indeed we have an exact sequence
\begin{equation}
1\rightarrow K^*_+\rightarrow K^*\xrightarrow{\ s\ } \{\pm 1\}\times \{\pm 1\}\rightarrow 1 \,,\  \mbox{ where } \ s(\alpha )=(\operatorname{sgn}(\alpha ),\operatorname{sgn}(\sigma (\alpha )))\,;  \end{equation}
since $-1\in\O^*$, the subgroup $s(\O^*)$ is either $\{\pm 1\} $ or $\{\pm 1\}\times \{\pm 1\}$, and the latter occurs if and only if there is a positive unit $\alpha $ with $\sigma (\alpha )<0$, equivalently with $\operatorname{Nm}(\alpha )=-1 $.

\bigskip	
\section{The key lemma}
We work over the complex numbers. 
\begin{lem}
Let $X,S$ be smooth projective varieties, and $\pi :X\rightarrow S$ a  double covering. Let $(E_i)_{i\in R}$ be the connected components of the branch locus of $\pi $. Let $(e_i)_{i\in R}$ be  the canonical basis of $(\Z/2)^R$, $\mathbf{e}:=\sum e_i$,
and let $\varphi : (\Z/2)^R\rightarrow \Pic(S)\otimes (\Z/2)$ be the homomorphism such that $\varphi (e_i)=[E_i]$. Then
$\Pic(X)[2]$ is isomorphic to $\Ker \varphi /\lr{\mathbf{e}}$.
\end{lem}
The inverse isomorphism is given as follows:  let $I\subset R$ such that $\sum\limits_{i\in I}e_i\in\Ker \varphi $, so that   $\sum\limits_{i\in I}E_i\equiv 2L $ for some $L\in\Pic(S)$. Put $F_i:=\pi ^{-1}(E_i)$; then $\sum\limits_{i\in I}F_i-\pi ^*L$ is the element of $\Pic(X)[2]$ corresponding to  $\sum\limits_{i\in I}e_i$. 

\medskip	
The lemma can be applied in two ways: information on the branch locus can give information on $\Pic(X)[2]$, or the other way around.  We will explain in the next section how even sets of nodes give an example of the latter situation. Let us give first an example of the former.

\medskip	
\noindent\emph{Example : Campedelli surfaces}.

In\cite{Cam}, Campedelli claimed to construct a surface of general type with $p_g=0$ and  $K^2=1$ by taking a double covering of $\P^2$ branched along a curve $C$ of degree 10 with one ordinary quadruple point $q$, 5 triple points $p_1,\ldots ,p_5$ of type $(3,3)$\footnote{ This means that they become  ordinary triple points after one blow up.}, and no other singularity.  

Indeed, assume that such a curve exists. Let $b:S\rightarrow \P^2$ be the blowing up of $q,p_1,\ldots ,p_5$, then  of the intersection point of the exceptional divisor above each $p_i$ with the strict transform of $C$. 
The fiber $b^{-1}(p_i)$ is the union of the new exceptional divisor $E'_i$ and the strict transform $\tilde{E}_i $ of the exceptional divisor above $p_i$. 
 The strict transforms $\tilde{C}$ of $C$ and the curves $\tilde{E}_i $ are smooth and do not meet. Let $F$ be the exceptional divisor above $q$, and let $L$ be a line in $\P^2$. We have  
$\ \tilde{C} \equiv 10\,b^*L-4F-\sum_{i=1}^5 (3\tilde{E}_i+6E'_i) \ $ in $\Pic(S)$.
It follows that the divisor $R:=\tilde{C}+\sum \tilde{E}_i$ is 
divisible by 2 in $\Pic(S)$, hence we can form the double covering $\pi :\hat{X}\rightarrow S$ branched along $R$. The $(-2)$-curves $\tilde{E}_i$ become $(-1)$-curves on $\hat{X}$, and we can contract them and get a minimal surface $X$. If the six points $q,p_1,\ldots ,p_5$ do not lie on a conic, it is easy to check that $p_g(X)=0$ and $K_X^2=1$.

Unfortunately Campedelli's construction of the curve $C$ was incorrect. A correct one appears in \cite{O-P}, and another one in \cite{W} (using some computer algebra). 
Werner takes for $C$ the union of a plane curve $B$ of degree 8 and a conic $Q$. The points $q$ and $p_5$ are  on $B\smallsetminus Q$; the points $p_1,\ldots ,p_4$ are on $Q$, and are tacnodes for $B$, with the same tangents for $B$ and $Q$. Let $\tilde{B} $ and $\tilde{Q} $ be the strict transforms of $B$ and $Q$ in $S$; we have $R=\tilde{B}+\tilde{Q}+\sum \tilde{E_i}   $. But $\tilde{Q}+\sum\limits_{i=1}^4\tilde{E}_i  \equiv 2b^*L$ is divisible by 2 in $\Pic(S)$, hence by the key lemma $\Pic(\hat{X})[2]=\Pic(X)[2]$ is nonzero.  The torsion subgroup of $ \Pic(X)$ is cyclic, of order $\leq 5$ \cite{Mi}; therefore it is $\Z/2$ or $\Z/4$. Werner shows that the linear system $\abs{3K_X}$ is base point free; by \cite{Mi} this implies $\operatorname{Tors}(\Pic(X))=\Z/2 $. 

The same method applies to the Oort-Peters construction. In that case  $\abs{3K_X}$ has a base point, so  $\operatorname{Tors}(\Pic(X))=\Z/4 $ \cite{W}.

\bigskip	
\section{Even sets of nodes}
I used  the key lemma to tackle a classical problem: what is the maximum number
$\mu (d)$ of nodes (=~ordinary double points) that a nodal surface of degree $d$ in $\P^3$ can have? It is classical that $\mu (3)=4$, the maximum being realized by the Cayley cubic surface $\sum\dfrac{1}{X_i}=0 $; and that $\mu (4)=16$, the maximum being realized by the Kummer surface. In \cite{S} Severi claimed the inequality   $\mu (5)\leq 31$, which implies $\mu (5)=31$ because Togliatti has constructed a quintic surface with 31 nodes \cite{T}; but the proof was insufficient.

 I proved this inequality in 1979 by applying, among others, the key lemma \cite{Be}. The same method, with some harder work, allowed Jaffe and Ruberman to prove $\mu (6)=65$ \cite{J-R}; here the maximum is achieved by a sextic surface constructed by Barth \cite{Ba}. The problem is still wide open for $d\geq 7$.

The basic idea is that a large set of nodes on a surface must contain some particular subsets, called even, which are easier to control.
Let $\Sigma $ be a projective surface, with  a finite set $\mathscr{N}$ of nodes and no other singularities. Let $b:S\rightarrow \Sigma $ be the minimal resolution of $\Sigma $; for $n\in\mathscr{N}$, $b^{-1}(n)$ is a \hbox{$(-2)$-curve} $E_n$. We say that a subset $R\subset \mathscr{N}$ is \emph{even}
if $\sum\limits_{n\in R}E_n$ is divisible by 2 in $\Pic(S)$; equivalently, there exists a double covering $X\rightarrow S$ branched along $\displaystyle\bigcup_{n\in R}E_n$ (this notion was introduced and thoroughly studied by Catanese in \cite{Ca}, see also \cite{C-C}). In this situation Riemann-Roch gives $\chi (\O_X)=2\chi (\O_S)-\dfrac{r}{4} $, with $r:=\operatorname{Card}(R) $; in particular, $r$ is divisible by 4. 

Let us now specialize to the case of a nodal quintic surface $\Sigma \subset\P^3$.
An inequality of Castelnuovo gives $r\geq 16$; and there are indeed some natural constructions of even sets with $r=16$, and also $r=20$. If $r\geq 24$, we get $\chi (\O_X)< \chi (\O_S)$, hence $q(X)\geq 1$; therefore $\Pic(X)[2]\neq 0$, so by the lemma our subset $R$ is the disjoint union of two even subsets, hence $r\geq 2\times 16=32$.

Suppose that $\Sigma $ has (at least) 32 nodes; the corresponding $(-2)$-curves $E_n$ define a map \break$(\Z/2)^{32}\rightarrow H^2(S,\Z/2)$. The image is totally isotropic in $H^2(S,\Z/2)$, hence of dimension $\leq \bigl[\dfrac{1}{2}b_2(S)\bigr]\allowbreak =26$. Therefore its kernel  has dimension $\geq 6$. 
This kernel is a vector subspace of $(\Z/2)^{32}$ whose elements have weight (= number of nonzero coordinates) $0,16,20$ or $32$. Some easy linear algebra over $\Z/2$ shows that this is impossible.

\bigskip	
\section{Proof of the key lemma}

Let $K_V$ be the function field of a variety $V$.
We start from the exact sequence
\[1\rightarrow \C^*\rightarrow K_{S}^*\xrightarrow{\ \operatorname{div} \ }\operatorname{Div}(S)\rightarrow \Pic(S) \rightarrow 0\]
which we split as
\[1\rightarrow \C^*\rightarrow K_{S}^*\rightarrow K_S^*/\C^*\rightarrow 1\quad ,\quad 1\rightarrow K_S^*/\C^*\rightarrow\operatorname{Div}(S)\rightarrow \Pic(S) \rightarrow 0\, .\]
The group $G=\Z/2$ acts on $X$ via the involution $\sigma $ of $X$ which swaps the two sheets of $\pi $.
We compare the second exact sequence with the $G$-invariants of the corresponding sequence for $X$:
\begin{equation}\label{diagram}
\begin{gathered}
\xymatrix{1 \ar[r] & K_S^*/\C^* \ar[r]\ar[d]^{\alpha } &\Div(S)  \ar[r]\ar[d]^{\beta } &\Pic(S) \ar[r]\ar[d]^{\gamma } & 0 \\
1 \ar[r] & (K_X^*/\C^*)^G \ar[r]&\Div(X)^G  \ar[r] &\Pic(X)^G\ar[r]  &H^1(G, K_X^*/\C^*) \,.}\end{gathered}\end{equation}
The essential points of the proof are the two following cohomological facts:

\smallskip	
\noindent\textbf{Fact 1}:  $H^1(G, K_X^*/\C^*)=0$: this follows from $H^1(G,K_X^*)=0$ (Hilbert's theorem 90) and $H^2(G,\C^*)=\allowbreak \C^*/\C^{*2}=0$.

\smallskip
\noindent\textbf{Fact 2}: $\operatorname{Coker}\alpha =\Z/2 $: this follows from 
 $H^1(G,\C^*)=\Ker \bigl(\C^*\xrightarrow{\ \times 2\ }\C^*\bigr)=\Z/2$ and the diagram
 \[\xymatrix{1 \ar[r] & \C^* \ar[r]\ar@{=}[d] &K_S^*  \ar[r]\ar@{=}[d] &K_S^*/\C^*\ar[r]\ar[d]^{\alpha } & 1 &\\
1 \ar[r] & \C^* \ar[r]&(K_X^*)^G  \ar[r] & (K_X^*/\C^*)^G \ar[r] & H^1(G,\C^*)\ar[r]&0\, .}\]

\smallskip	
Now we consider the exact sequence of cokernels in the diagram (2). Since $\Pic(S)[2]=(0)$,  $\gamma $ is injective.
$\operatorname{Div}(X)^G $ is the free abelian group with basis the divisors $C+\sigma  (C)=\pi ^*C$, for $C$ an irreducible curve on $S$ distinct from the $E_i$, and the curves $F_i:=\pi ^{-1}(E_i)$. Therefore $\operatorname{Coker}\beta  =(\Z/2)^{R}$, 
and we get an exact sequence
\[0\rightarrow \Coker \alpha =\Z/2 \xrightarrow{\ j\ } \Coker \beta =(\Z/2)^{R}\xrightarrow{\ k\ } \operatorname{Coker}\gamma \rightarrow 0 \,,\]with   $k(e_i)=[F_i]$ in $\Pic(X)^G\pmod{\im \gamma }$. 
Since $\sum E_i$ is divisible by 2 in $\Pic(S)$, $k$ maps  $\mathbf{e}:=\sum e_i$ to $0$, hence $j(1)=\mathbf{e}$. Thus we get an isomorphism $\operatorname{Coker}\gamma \cong  (\Z/2)^R/\lr{\mathbf{e}}$, with $\mathbf{e}:=\sum e_i$. 

\smallskip	
Now we apply  the snake lemma to the diagram
\[\xymatrix{0\ar[r]& \Pic(S)\ar[r]^<<<<{\gamma }\ar[d]^{\times 2}& \Pic(X)^G \ar[r]\ar[d]^{\times 2}&(\Z/2)^R/\lr{\mathbf{e}}  \ar[r]\ar[d]^{\times 2}& 0\,{}\\
0\ar[r]& \Pic(S)\ar[r]^<<<<{\gamma }& \Pic(X)^G \ar[r] &(\Z/2)^R/\lr{\mathbf{e}} \ar[r]& 0\,.
}\]We get an exact sequence $0\rightarrow \Pic(X)^G[2]\rightarrow  (\Z/2)^R/\lr{\mathbf{e}}\xrightarrow{\ \bar{\varphi }\ } \Pic(S)\otimes  (\Z/2)$, with $\bar{\varphi }(e_i)=[E_i]$. Observe that $\Pic(X)^G[2]=\Pic(X)[2]$: indeed 
for $L\in\Pic(X)[2]$ we have $\operatorname{Nm}L\in \Pic(S)[2]=(0)$, hence $\sigma ^*L\otimes L=\O_X$, so $\sigma  ^*L=L^{-1}=L$. Finally we get 
our isomorphism $\psi:\Pic(X)[2]\iso \Ker \varphi /\lr{\mathbf{e}}$. 

Let $\mathbf{f}=\sum\limits_{i\in I}e_i$ be an element of $\Ker \varphi $; then
$\sum\limits_{i\in I}E_i\equiv 2L$ in $\Pic(S)$, and $\sum\limits_{i\in I}F_i-\pi ^*L$ is an element of $\Pic(X)[2]$ which maps to $\mathbf{f}$ under $\psi$. \qed
  
\bigskip	
\noindent\emph{Remark}$.-$ With no hypothesis on $\Pic(S)[2]$, the proof gives a (split) exact sequence\footnote{ We assume here $R\neq \varnothing$; when $\pi $ is \'etale $\Ker\pi ^*$ has order $2$.}\[0\rightarrow \Pic(S)[2]\xrightarrow {\ \pi ^*\ }\Pic(X)^G[2]\rightarrow \Ker \varphi /(1,\ldots ,1)\rightarrow 0\]and an isomorphism  $\Pic(X)^G[2]= \Ker \bigl(\Pic(X)[2]\xrightarrow{\ \operatorname{Nm} \ }\Pic(S)[2]\bigr)$. 
In other words, the line bundles $L$ on $X$ with $L^2=\O_X$ and $\operatorname{Nm}L=\O_S $ are of the form $\pi ^*M(-\sum\limits_{i\in I} F_i)$ for some subset $I\subset R$ and line bundle $M$ on $S$ such that $M^{2}=\O_S(\sum\limits_{i\in I}E_i)$;    this presentation is unique up to replacing $(I,M)$ by $(R\smallsetminus I, (\mathscr{L}\otimes M)^{-1})$, with $\mathscr{L}=\det(\pi _*\O_X)$.

\smallskip	
This applies in particular when $\pi $ is a double covering of projective curves; we get a description of $P[2]$, where $P:=\Ker \bigl(JX\xrightarrow{\ \operatorname{Nm} \ }JS\bigr)$ is the Prym variety of $(X,S)$. Note that in that case $\Pic(S)\otimes (\Z/2)=\Z/2$, so the the subsets $I\subset R$ which appear in that description are those with  $\# I$ even.

\bigskip
\section{Proof of Gauss' theorem}

Let $d$ be a square-free integer, and let $\O$ be the ring of integers of $K:=\Q(\sqrt{d})$. The   proof of the key lemma can be easily adapted to the case $S=\operatorname{Spec}\Z $, $X=\operatorname{Spec}\O$: we replace the diagram (2) by 
\[\xymatrix{1 \ar[r] & \Q^*_+\ar[r]\ar[d]^{\alpha } &\Div(\Z)  \ar[r]\ar[d]^{\beta } &\Pic_+(\Z)=0 \ar[r]\ar[d]^{\gamma } & 0 \\
1 \ar[r] & (K_+^*/\O_+^*)^G \ar[r]&\Div(\O)^G  \ar[r] &\Pic_+(\O)^G\ar[r]  &H^1(G, K_+^*/\O_+^*) \,.}\]
We just have to check that the two cohomological facts used in the previous proof still hold.
We recall that if $M$ is a $G$-module, that is, an abelian group with an involution $\sigma $, we have 

\centerline{$H^{1}(G,M)=\Ker(1+\sigma )/\im(1-\sigma )\ $ and $\ H^{2}(\sigma,M)=\Ker(1-\sigma )/\im(1+\sigma )$. }

\smallskip	
\noindent\textbf{Fact 1}:  $H^1(G,K^*_+/\O^*_+)=0$. 

We use the exact sequence $H^1(G,K^*_+)\rightarrow H^1(G,K^*_+/\O^*_+)\rightarrow H^2(G,\O^*_+)\rightarrow H^2(G,K^*_+)$.
If we make $G $ act on $\{\pm 1\}\times \{\pm 1\}$ by swapping the factors, the  sequence $(1)$ of \S\,2  is an exact sequence of $G$-modules. Taking invariants gives an exact sequence $\Q^*\xrightarrow{\ \operatorname{sgn}\ } \{\pm 1\}\rightarrow H^1(G,K^*_+)\rightarrow H^1(G,K^*) $; since $H^1(G,K^*)=0$ by Hilbert's Theorem 90, we get $H^1(G,K^*_+)=0$.

We have $H^2(G,\O^*_+)=(\O^*_+)^G/\operatorname{Nm}\O^*_+= (\O^*_+)^G$. This group is trivial if $d>0$, and equal to  $\{\pm 1\} $ if $d<0$.  In the latter case the map $H^2(G,\O^*)\rightarrow H^2(G,K^*)$ is the composite $\{\pm 1\}\hookrightarrow \Q^*\rightarrow \Q^*/\operatorname{Nm}(K^*)   $. Since $\operatorname{Nm} (K^*)\subset \Q_+$ this map is injective. So in both cases we find  $H^1(G,K^*_+/\O^*_+)=0$.

\medskip	
\noindent\textbf{Fact 2}: $H^1(G,\O^*_+) =\Z/2$. 

If $d>0$, we have $\O^*_+=\Z$, with $G $ acting by changing sign; thus $H^1(G,\O^*_+)=\Z/2$. If $d<-3$, $\O^*_+=\O^*=\{\pm 1\} $, hence the result. For $d=-1$ we observe that $(1-\sigma )(i)= -1$, and for $d=-3$ $(1-\sigma )(\rho )=\rho ^2$, so $H^1(G,\O^*) =\O^*/\im(1-\sigma )=\{\pm 1\} $.

\medskip	
Thus as before we get an exact sequence
\[0\rightarrow \Coker \alpha =\Z/2 \xrightarrow{\ j\ } \Coker \beta =(\Z/2)^{R}\xrightarrow{\ k\ } \operatorname{Coker}\gamma =\Pic_+(\O)^G\rightarrow 0 \,,\]where $R\subset \operatorname{Spec}\Z $ is the set of ramified primes. As  above we have $\Pic_+(\O)^G=\Pic_+(\O)[2]=\operatorname{Cl}_+(K) $, hence Gauss' theorem -- in fact we get as above an explicit isomorphism $(\Z/2)^R/\lr{\mathbf{e}}\iso \operatorname{Cl}_+(K) $.\qed


\bigskip	
\begin{thebibliography}{X-X}

\bibitem[Ba]{Ba} W. Barth\,: \textsl{Two projective surfaces with many nodes, admitting the symmetries of the icosahedron}. 
J. Algebraic Geom. \textbf{5}  (1996), no. 1, 173-186.

\bibitem[Be]{Be} A. Beauville\,: \textsl{Sur le nombre maximum de points doubles d'une surface dans} $\P^3$ ($\mu(5)=31$).  Journ\'ees de G\'eom\'etrie Alg\'ebrique d'Angers, pp. 207-215, Sijthoff \& Noordhoff (1980).   


\bibitem[C]{Cam} L. Campedelli\,: \textsl{Sopra alcuni piani doppi notevoli con curva di diramazioni de1 decimo
ordine}. Atti Acad. Naz. Lincei \textbf{15} , 536-542 (1932).
 
 \bibitem[C-C]{C-C} G. Casnati, F. Catanese\,: \textsl{Even sets of nodes are bundle symmetric}. J. Differential Geom. \textbf{47}  (1997), no. 2, 237-256.
 
\bibitem[Ca]{Ca} F. Catanese\,: \textsl{Babbage's conjecture, contact of surfaces, symmetric determinantal varieties and applications}. Invent. Math. \textbf{63}  (1981), no. 3, 433-465.

\bibitem[J-R]{J-R} D. Jaffe, D. Ruberman\,: \textsl{A sextic surface cannot have $66$ nodes}. 
J. Algebraic Geom. \textbf{6}  (1997), no. 1, 151-168. 


\bibitem[Mi]{Mi} Y. Miyaoka\,: \textsl{Tricanonical maps of numerical Godeaux surfaces}. Invent. Math. \textbf{34} (1976), no. 2, 99-111.

\bibitem[O-P]{O-P} F. Oort, C. Peters\,: \textsl{A Campedelli surface with torsion group} $\Z/2$. 
Nederl. Akad. Wetensch. Indag. Math. \textbf{43}  (1981), no. 4, 399-407. 

\bibitem[S]{S} F. Severi\,: \textsl{Sul massimo numero di nodi di una superficie di dato ordine dello spazio ordinario o di una forma di un iperspazio}.  Ann. Mat. Pura Appl. (4) \textbf{25}  (1946), 1-41. 

\bibitem[T]{T} E. Togliatti\,: \textsl{Una notevole superficie del $5^{\rm o}$ ordine con soli punti doppi isolati}. 
Vierteljschr. Naturforsch. Ges. Z\"urich \textbf{85}  (1940), Beiblatt (Festschrift Rudolf Fueter), 127-132. 


\bibitem[W]{W} C. Werner\,: \textsl{A surface of general type with} $p_g=q=0$, $K^2=1$. Manuscripta Math. \textbf{84}  (1994), no. 3-4, 327-341.
\end{thebibliography}
\end{document}